\documentclass[12pt,a4paper]{amsart}
\usepackage[utf8]{inputenc}	
\usepackage[T1]{fontenc}
\usepackage[in,headings]{fullpage}
\usepackage{amsbsy, amscd, amsfonts, amsmath, amsrefs, amssymb, amsthm}
\usepackage{graphicx}
\usepackage{float}
\usepackage{color}   
\usepackage[colorlinks=true,linkcolor=blue,citecolor=blue,urlcolor=blue]{hyperref}

\newtheorem{theorem}{Theorem}

\newtheorem{proposition}[theorem]{Proposition}
\newtheorem{lemma}[theorem]{Lemma}

\theoremstyle{definition}
\newtheorem{definition}[theorem]{Definition}
\newtheorem{remark}[theorem]{Remark}

\theoremstyle{remark}

\newcommand{\C}{\mathbf{C}}
\newcommand{\Z}{\mathbf{Z}}
\newcommand{\Q}{\mathbf{Q}}
\newcommand{\R}{\mathbf{R}}
\newcommand{\N}{\mathbf{N}}

\newcommand{\Hip}{\mathbb{H}}
\renewcommand{\Re}{\mathop{\mathrm{Re}}\nolimits}
\renewcommand{\Im}{\mathop{\mathrm{Im}}\nolimits}
\newcommand{\Rzeta}{\mathop{\mathcal R }\nolimits}

\newfont{\cmbsy}{cmbsy10}
\newfont{\cmmib}{cmmib10}
\newcommand{\Orden}{\mathop{\hbox{\cmbsy O}}\nolimits}

\DeclareMathOperator*{\Res}{Res}

\begin{document}

\title[Levinson Functions]
{Levinson Functions.}
\author[Arias de Reyna]{J. Arias de Reyna}
\address{%
Universidad de Sevilla \\ 
Facultad de Matem\'aticas \\ 
c/Tarfia, sn \\ 
41012-Sevilla \\ 
Spain.} 



\email{arias@us.es, ariasdereyna1947@gmail.com}


\begin{abstract}
Starting from some of Norman Levinson's results, we construct interesting examples of functions $f(s)$ such that for $s=\frac12+it$, we have $Z(t)=2\Re\{\pi^{-\frac{s}{2}}\Gamma(s/2)f(s)\}$. For example one such function is 
\[
\begin{aligned}
\Rzeta_{-3}(s)=\frac12&\int_{0\swarrow1}\frac{x^{-s}e^{3\pi ix^2}}{e^{\pi i x}-e^{-\pi i x}}\,dx\\&+\frac{1}{2\sqrt{3}}\int_{0\swarrow1}\frac{x^{-s}e^{\frac{\pi i}{3}x^2}}{e^{\pi i x}-e^{-\pi i x}}\Bigl(e^{\frac{\pi i}{2}}+2e^{-\frac{\pi i}{6}}\cos(\tfrac{2\pi x}{3})\Bigr)\,dx.
\end{aligned}
\]
\end{abstract}

\maketitle

\section{Introduction} 
The functional equation for the zeta function is  $h(s)\zeta(s)=h(1-s)\zeta(1-s)$ with $h(s)=\pi^{-s/2}\Gamma(s/2)$. In the critical line, we have 
\begin{equation}
h(\tfrac12+it)=e^{i\vartheta(t)}|h(\tfrac12+it)|.
\end{equation}
where $\vartheta(t)$ is a real and real analytic function defined in $\R$, which can be extended as an analytic function to the complex plane with two cuts along $[1,\infty)$ and $(-\infty, 0]$.  The functional equation implies that $-\vartheta(t)$ is also the phase of zeta in the critical line. That is, we have (see \cite{T}*{4.17})
\begin{equation}
\zeta(\tfrac12+it)=e^{-i\vartheta(t)} Z(t),\qquad t\in\R,
\end{equation}
where $Z(t)$ is a real and real analytic function in $\R$, the function $Z(t)$ extends to an analytic function in the same region as $\vartheta$.

The auxiliary function of Riemann $\Rzeta(s)$ satisfies the equation 
\[Z(t)=2\Re\{\pi^{-s/2}\Gamma(s/2)\Rzeta(s)\},\qquad s=\tfrac12+it,\quad t\in\R.\]
So, the ordinates of the zeros of $\zeta(s)$ on the critical line are determined by $\Rzeta(s)$. There are many functions that share this property with $\Rzeta(s)$ and we call them Levinson functions.
In Section \ref{S:dos} we give the definition, give several interesting examples, and define the related concept of Levinson pair.

Levinson \cite{L2} constructed interesting examples of Levinson pairs. In Section \ref{S:tres} we construct from these Levinson pairs a family $\Rzeta_\tau(s)$ of Levinson functions one for each $\tau\ne0$ in the closed upper half-plane $\Hip$. This extends the example of Riemann, since $\Rzeta_{-1}(s)=\Rzeta(s)$. 

The functions $\Rzeta_\tau(s)$ for $\tau\in\Q^*$ can be given in a simplified form. We dedicate several sections to get this simplified form. This included an application of Gauss sums. The results are nice expressions for $Z(t)$, for example,
\[Z(t)=2\Re\bigl\{e^{i\vartheta(t)} \Rzeta_{-3}(\tfrac12+it)\Bigr\},\]
with 
\[
\begin{aligned}
\Rzeta_{-3}(s)=\frac12&\int_{0\swarrow1}\frac{x^{-s}e^{3\pi ix^2}}{e^{\pi i x}-e^{-\pi i x}}\,dx\\&+\frac{1}{2\sqrt{3}}\int_{0\swarrow1}\frac{x^{-s}e^{\frac{\pi i}{3}x^2}}{e^{\pi i x}-e^{-\pi i x}}\Bigl(e^{\frac{\pi i}{2}}+2e^{-\frac{\pi i}{6}}\cos(\tfrac{2\pi x}{3})\Bigr)\,dx.
\end{aligned}
\]
From these formulas, we may  get a version of the approximate functional equation for $\zeta(s)$ \cite{T}*{Thm. 4.13} in which the error term is substituted by an asymptotic expansion of Riemann-Siegel type.

In the end, we include the x-rays of some of these functions. It will be interesting to study how the zeros of $\Rzeta_\tau(s)$ move when $\tau$ varies.

\subsection{Notations}

We consider meromorphic functions $f(z)$ defined in the complex plane. Given such a function, it is well known that the conjugate function $f^*(z):=\overline{f(\overline z)}$ is also a meromorphic function. 
We have $f=f^*$ if and only if $f$ takes real values on the real line.

For integers $a$ and $b$, we use the notation $a\perp b$ to indicate that $\gcd(a,b)=1$.

\section{Definition of Levinson functions}\label{S:dos}

\begin{definition}
We say that a meromorphic function $f$ defined in $\C$ is a Levinson function if we have for $t\in\R$
\begin{equation}
Z(t)=2\Re\bigl\{e^{i\vartheta(t)} f(\tfrac12+it)\Bigr\}.
\end{equation}
\end{definition}

Notice that this is equivalent to 
\[e^{i\vartheta(t)}\zeta(\tfrac12+it)=e^{i\vartheta(t)} f(\tfrac12+it)+e^{-i\vartheta(t)} \overline{f}(\tfrac12-it).\]
And this is equivalent to saying that for $s$ in the critical line
\begin{equation}\label{E:2def}
h(s)\zeta(s)=h(s)f(s)+h(1-s)f^*(1-s).
\end{equation}
Then \eqref{E:2def} is true for all $s$. And a meromorphic function $f$ is a Levinson function  if and only if it satisfies  \eqref{E:2def}. 

Since $\zeta(s)$ is real for real $s$, we have $\zeta^*=\zeta$, and then the functional equation shows that $\zeta(s)/2$ is a Levinson function. But we have more interesting examples. 

\subsection{Example of Levinson}
In his proof that more than $\frac13$ of the zeros of zeta are on the critical line, Levinson \cite{L} proved that 
\begin{equation}
L(s)=\zeta(s)+\frac{\zeta'(s)}{\frac{h'(1-s)}{h(1-s)}+\frac{h'(s)}{h(s)}},
\end{equation}
is a Levinson function. Conrey \cites{C1, C2, C3} gives other related examples. 

In \cite{A62} we show, assuming RH, that the phase of $L(\frac12+it)$ is $2\pi-E(t)$ where $E(t)$ is a real analytic version of $S(t)$ so that $|\pi S(t)-E(t)|\le\pi/2$.  The difference of $\pi S(t)-E(t)$ counting the failures of RH  when RH is not assumed.

\subsection{Examples of Riemann}
(a) The first interesting Levinson function was defined by Riemann in his paper \cite{R}. He proved the following
\begin{equation}
\pi^{-s/2}\Gamma(s/2)\zeta(s)=\frac{1}{s(s-1)}+\int_1^{+\infty}(x^{s/2}+x^{(1-s)/2})\frac{\theta(x)-1}{2}\frac{dx}{x},
\end{equation}
where
\begin{equation}
\theta(x)=\sum_{n\in\Z}e^{-\pi n^2 x}.
\end{equation}
This is equivalent to saying that the function $\Rzeta_i(s)$ defined by 
\begin{equation}\label{E:1stRiemann}
\pi^{-s/2}\Gamma(s/2)\Rzeta_i(s)=-\frac1s+\int_1^{+\infty}x^{s/2}\;\frac{\theta(x)-1}{2}\frac{dx}{x},
\end{equation}
with ${\Rzeta_i}^*=\Rzeta_i$, is a Levinson function. (The notation used for this function will be explained later). 

(b) Riemann found another Levinson function. He did not publish anything about it, but Siegel \cite{Siegel} was able to extract it from the papers Riemann left behind after his death. The function
\begin{equation}
\Rzeta(s)=\int_{0\swarrow1}\frac{x^{-s} e^{\pi i x^2}}{e^{\pi i x}-e^{-\pi i x}}\,dx,
\end{equation}
is a Levinson function. Here, the integration is, for example, along the line parametrised by $u\mapsto\frac12+ u e^{-3\pi i/4}$.

\section{Levinson pairs}
Before giving this generalisation, it is convenient to define Levinson pairs. 
We say that a pair $(f,g)$ of meromorphic functions defined on $\C$ is a \emph{Levinson pair} if they satisfy the relation
\begin{equation}\label{D:LevPair}
h(s)\zeta(s)=h(s)f(s)+h(1-s)g(1-s).
\end{equation}
To say that $f$ is a Levinson function is the same as say that $(f,f^*)$ is a Levinson pair. The next two propositions shows that being a Levinson pair, by itself, is not a particularly interesting property.
\begin{proposition}\label{P:uno}
Let $f$ be a meromorphic function defined on the complex plane. There exists a meromorphic function $g$ such that $(f,g)$ is a Levinson pair. 
\end{proposition}
\begin{proof}
Such a function $g$ will satisfy \eqref{D:LevPair}, so that 
\begin{equation}\label{E:par1}
h(1-s)\zeta(1-s)=h(1-s)f(1-s)+h(s)g(s).
\end{equation}
And there is one and only one meromorphic function $g$ satisfying this 
\begin{equation}
g(s):=\chi(s)\bigl(\zeta(1-s)-f(1-s)\bigr), 
\end{equation}
where $\chi(s)=h(1-s)/h(s)$ is the meromorphic function appearing in the functional equation 
$\zeta(s)=\chi(s)\zeta(1-s)$.
\end{proof}

\begin{proposition}\label{P:properties}
If $(f,g)$ is a Levinson pair so it is $(g,f)$, and $(f^*,g^*)$.
\end{proposition}
\begin{proof}
We have just seen that if $(f,g)$ is a Levinson pair, then we have \eqref{E:par1}. Therefore, by the functional equation of $\zeta(s)$ we get 
\[h(s)\zeta(s)=h(s)g(s)+h(1-s)f(1-s).\]
So, $(g,f)$ is a Levinson pair. 

Notice that a meromorphic function $f$ is real on the real axis if and only if $ f^*=f$. Therefore, from \eqref{D:LevPair} we easily get that 
\[h(s)\zeta(s)=h(s)f^*(s)+h(1-s)g^*(1-s).\]
(Notice that $(fg)^*=f^* g^*$ and $(f+g)^*=f^*+g^*$).
\end{proof}

\begin{proposition}\label{P:dos}
Let $(f,g)$ be a Levinson pair, then $(f+g^*)/2$ is a Levinson function.
\end{proposition}
\begin{proof}
Notice that the set of Levinson pairs $(f,g)$ is a convex set with the usual sum and scalar product. So, by Proposition  \ref{P:properties} the pair \[((f+g^*)/2, (g+f^*)/2)\] is a Levinson pair. But it is clear that if $v=\frac12(f+g^*)$, we have $ v^*=\frac12(f^*+g)$. Therefore, $(v,v^*)$ is a Levinson pair, which means that $v$ is a Levinson function.
\end{proof}

Propositions \ref{P:uno} and \ref{P:dos} give us cheap but uninteresting Levinson functions.

\section{Levinson generalisation of Riemann's examples}\label{S:tres}

(d) Our main examples are given by Levinson in \cite{L2}, and it is a generalisation of 
the examples of Riemann. First following  Levinson we consider its function
\begin{equation}
f(s,\tau)=\frac{1}{2\pi i}\int_{L} e^{-\pi i \tau x^2}\Bigl(\sum_{n\in\Z}\frac{x e^{\pi i \tau n^2}}{x^2-n^2}\Bigr)x^{-s}\,dx.
\end{equation}
For the convergence of the sum in the integrand, we need to assume $\Im(\tau)\ge0$. 
Given $\tau$ satisfying this condition, we must integrate along a line $L$ passing through the point $x=\frac12$ and a slope such that $\Re(\pi i\tau x^2)>0$ for $|x|$ large. This can be achieved, for example, taking $L$ with parametrisation $x=\frac12+u e^{i\theta}$ with $\theta=-\frac{\pi}{4}-\frac{\alpha}{2}$, $0<\alpha<\pi$ being the argument of $\tau$. We always take $L$ oriented so that the imaginary part of $x$ decreases along the line $L$. 

The function  $f(s,\tau)$ is an entire function of $s$ for each fixed value $\tau\in\overline{\Hip}^*$ where $\Hip$ is the upper  hyperplane, and $\overline{\Hip}^*$ its closure excluding the value $\tau=0$. 
In particular, we have $f(s,-1)=\Rzeta(s)$. For $s$ fixed, $f(s,\tau)$ is holomorphic in the open set $\Im\tau>0$ and continuous in $\Hip^*$.  All this is proved in \cite{L2}\footnote{Our parameter $\tau$ does not match Levinson's $\alpha$, they are related by $\tau=i\alpha$.}.

\begin{lemma}\label{Lf*}
For any $\tau\in\overline \Hip^*$ and $s\in \C$ we have
\begin{equation}
f^*(s,\tau)=f(s,-\overline\tau).
\end{equation}
\end{lemma}

\begin{proof}
Consider first the case $\Im\tau>0$ and $\Re s<0$,  we have
\[f(\overline s,\tau)=\frac{1}{2\pi i}\int_{0\downarrow1}e^{-\pi i \tau x^2}\Bigl(\sum_{n\in\Z}\frac{x e^{\pi i \tau n^2}}{x^2-n^2}\Bigr)x^{-\overline s}\,dx.\]
We may take the path to be of direction $-i$, since then 
\[e^{-\pi i \tau x^2}=e^{-\pi i \tau (\frac12-ui)^2}=e^{-\frac{\pi i}{4}\tau-\pi\tau u+\pi i\tau u^2},\] 
and the factor $e^{-\pi\Im(\tau) u^2}$ makes the integral convergent. 

Since $\Re s<0$ we may move the line of integration to pass through the point $x=0$. Then with $x=-iu$
\[f(\overline s,\tau)=\frac{1}{2\pi i}\int_{-\infty}^\infty e^{\pi i \tau u^2}\sum_{n\in\Z}\frac{-iu e^{\pi i\tau n^2}}{-u^2-n^2}e^{-\overline s\log(-iu)}(-i\,du).\]
Taking complex conjugate of this expression 
\[f^*(s,\tau)=-\frac{1}{2\pi i}\int_{-\infty}^\infty e^{-\pi i \overline\tau u^2}\sum_{n\in\Z}\frac{iu e^{-\pi i\overline\tau n^2}}{-u^2-n^2}e^{-s\log(iu)}(i\,du).\]
Now take $x=iu$
\[f^*(s,\tau)=-\frac{1}{2\pi i}\int_{-i\infty}^{i\infty} e^{\pi i \overline\tau x^2}\sum_{n\in\Z}\frac{x e^{-\pi i\overline\tau n^2}}{x^2-n^2}x^{-s}\,dx=
f(s,-\overline\tau).\]
The equality then extends to all $s$ since both sides are entire functions of $s$. 
Then by continuity extends also to al $\tau\in\overline \Hip^*$.
\end{proof}

The main result of Levinson is that $(f(s,\tau), f^*(s,1/\overline{\tau}))$ is a Levinson pair
\begin{equation}\label{E:levmain}
h(s)\zeta(s)=h(s)f(s,\tau)+h(1-s)f^*(1-s,1/\overline\tau).
\end{equation}
By Proposition \ref{P:dos} we obtain that  
\begin{equation}\label{RAF}
\Rzeta_\tau(s):=\Rzeta(s,\tau):=\frac{f(s,\tau)+f(s,1/\overline\tau)}{2},
\end{equation}
is a Levinson function. 

\begin{definition}
For any $\tau\in\overline{\Hip}^*$,  the Levinson function defined in \eqref{RAF} will be called a Riemann auxiliary function.
\end{definition}

For $\tau\in\overline{\Hip}^*$ with  $|\tau|=1$ we have
$\tau=1/\overline{\tau}$ and it follows that in these cases $\Rzeta(s,\tau)=f(s,\tau)$ is a Levinson function.

In particular,  
\begin{equation}
\Rzeta(s)=\Rzeta(s,-1),\qquad  \Rzeta_i(s)=\Rzeta(s,i),
\end{equation} 
where $\Rzeta_i(s)$ is the first example of Riemann defined in \eqref{E:1stRiemann}. The first equality follows directly from the definitions, and the second follows from the results on the section about Lavrik's expression.

When $\tau\in\Q^*$ we are going to obtain simplified expressions of $\Rzeta_\tau(s)$, but we need some prior preparations.

\begin{proposition}
For $\tau\in\overline{\Hip}^*$ we have 
\begin{equation}
\Rzeta_{1/\overline\tau}(s)=\Rzeta_\tau(s),\quad \Rzeta_\tau^*(s)=\Rzeta_{-1/\tau}(s).\end{equation}
\end{proposition}
\begin{proof}
Direct from the Definition and Lemma \ref{Lf*}.
\end{proof}

\section{A Mittag-Leffler expansion}

\begin{lemma}\label{L:mittag}
Let $\alpha$ be a real number with $|\alpha|\le 1$, then
\begin{equation}\label{E:mit}
\lim_{N\to\infty}\frac{1}{2\pi i}\sum_{n=-N}^{n=N}(-1)^n\frac{z e^{\pi i n\alpha}}{z^2-n^2}=\frac12\frac{e^{\pi i \alpha z}+e^{-\pi i \alpha z}}{e^{\pi i z}-e^{-\pi i z}},
\end{equation}
\end{lemma}
\begin{proof}
For $\alpha=1$ and $\alpha=-1$ this is the well-known Mittag-Leffler expansion of $\cot\pi z$. When $|\alpha|<1$ let $f(z)=\frac12\frac{e^{\pi i \alpha z}+e^{-\pi i \alpha z}}{e^{\pi i z}-e^{-\pi i z}}$ and, given a natural number $N$, consider the integral 
\[I_N:=\frac{1}{2\pi i}\int_{\Gamma_N}\frac{f(\xi)}{\xi-z}\,d\xi,\]
where $\Gamma_N$ is the positively oriented contour of the rectangle $[-N-\frac12, N+\frac12]\times[-N,N]$. This integral is equal to the sum of residues at the poles in the rectangle
\begin{align*}
I_N&=\Res_{\xi=z}\frac{f(\xi)}{\xi-z}+\sum_{n=-N}^{k=N}\Res_{\xi=n}\frac{f(\xi)}{\xi-z}\\&=f(z)-\frac12\sum_{k=-N}^N\frac{e^{\pi i \alpha n}+e^{-\pi i \alpha n}}{2\pi i (-1)^n}\frac{1}{(z-n)}\\
&=f(z)-\frac{1}{4\pi i}\sum_{n=-N}^{n=N}(-1)^n\Bigl(\frac{e^{\pi i \alpha n}}{z-n}+
\frac{e^{\pi i \alpha n}}{z+n}\Bigr)\\
&=f(z)-\frac{1}{2\pi i}\sum_{n=-N}^N(-1)^n\frac{z e^{\pi i \alpha n}}{z^2-n^2}.
\end{align*}
So, our objective is to prove that $\lim_{N\to\infty}I_N=0$. We bound the integrals along the segments of the rectangle, and we may assume that $|z|\le N/2$
\[\Bigl|\int_{N+1/2-iN}^{N+1/2+iN}\frac{f(\xi)}{\xi-z}\,d\xi\Bigr|\le \frac{1}{N/2} \int_{-N}^N \frac{\cosh(\pi\alpha y)}{|(-1)^N ie^{-\pi y}+(-1)^N i e^{\pi y}|}\,dy=\Orden(\tfrac1N). \]
A similar bound is valid on the left vertical side. On the horizontal side of the rectangle, we have
\[\Bigl|\int_{-N-1/2-iN}^{N+1/2-iN}\frac{f(\xi)}{\xi-z}\,d\xi\Bigr|\le
\frac{1}{N/2}\int_{-N-1/2}^{N+1/2}\frac{\cosh(\pi\alpha N)}{e^{\pi N}-1}\,dy=\Orden(e^{-\pi(1-|\alpha|)N}).\]
On the other horizontal side, the situation is similar.
\end{proof}

\begin{remark}
With a simple computation, we can see that equation \eqref{E:mit} is equivalent to 
the Mittag-Leffler expansion 
\begin{equation}
\frac{1}{z}+\sum_{n\in\Z^*}(-1)^n \cos(\pi\alpha n)\Bigl(\frac{1}{z+n}-\frac{1}{n}\Bigr)=\pi\frac{\cos(\pi\alpha z)}{\sin\pi z},\qquad -1\le\alpha\le 1.
\end{equation}
\end{remark}

\section*{Mordell integral}
According to Siegel \cite{Siegel} Riemann considered the integral
\begin{equation}
\Phi(z,\tau):=\int_{0\uparrow1}\frac{e^{-\pi i\tau x^2+2\pi i z x}}{e^{\pi i x}-e^{-\pi i x}}\,dx,\qquad z\in\C,\quad \tau\in \Hip.
\end{equation}
Later, this integral was studied by Kronecker \cites{K1, K2}, Mordell \cite{M} and Zwegers \cite{Z}\footnote{Zwegers studies the Mordell integral $h(z,\tau)$, that he defines \[h(z,\tau)=\int_R\frac{e^{\pi i\tau x^2-2\pi zx}}{\cosh\pi x}\,dx.\] Our function is in Riemann's form. They are related by \[\Phi(z,\tau)=\tfrac12e^{-\frac{\pi i }{4}\tau+\pi i z}h(z-\tau/2;\tau).\]}. We have 
\begin{proposition}\label{P:Phi}
The function $\Phi(z,\tau)$ have the following properties
\begin{itemize}
\item[(1)] $\displaystyle{\Phi(z+1,\tau)-\Phi(z,\tau)=\frac{i}{\sqrt{-i\tau}}e^{\frac{\pi i}{\tau}(z+\frac12)^2}}$. 
\item[(2)]  $\displaystyle{\Phi(z+\tau,\tau)+e^{2\pi i z+\pi i\tau}\Phi(z,\tau)=e^{2\pi i z+\pi i\tau}.}$
\item[(3)] $z\mapsto \Phi(z,\tau)$ is the unique holomorphic function that satisfies (1) and (2) above.
\item[(4)] $\displaystyle{ \Phi\Bigl(\frac{z}{\tau}-\frac12,-\frac{1}{\tau}\Bigr)= -i\sqrt{-i \tau}\; e^{-\pi i z^2/\tau}\Phi\Bigl(z+\frac12,\tau\Bigr).}$
\end{itemize}
\end{proposition}

Applying Cauchy's Theorem to rotate the integration line, the  function $\Phi(z,\tau)$ with a fixed $z$, extends analytically to $-\pi<\arg\tau<2\pi$. Therefore, we may speak of the values when $\tau\in\R$.  Riemann was able to compute exact values for $\Phi(z,\tau)$ when $\tau$ is a rational number. These computations follow from  Proposition \ref{P:Phi}.
\begin{proposition}
Let $a$ and $b\in \N$, then 
\begin{equation}
\begin{aligned}[t]\label{E:integral1}
&(1-(-1)^{b(a+1)}e^{-2\pi i bz})\Phi(z,\tfrac{a}{b})=\\&\sum_{n=0}^{b-1}(-1)^ne^{-2\pi i n z}e^{-\pi i n^2 \frac{a}{b}}+(-1)^{b(a+1)}e^{-2\pi i bz}\frac{e^{3\pi i/4}}{\sqrt{a/b}}\sum_{m=0}^{a-1}e^{\frac{\pi i b}{a}(z+m+\frac12)^2},
\end{aligned}\end{equation}
\begin{equation}\label{E:integral2}
\begin{aligned}[t]
&(1-(-1)^{b(a+1)}e^{-2\pi i bz})\Phi(z,-\tfrac{a}{b})=\\&\sum_{n=0}^{b-1}(-1)^ne^{-2\pi i n z}e^{\pi i n^2 \frac{a}{b}}+(-1)^{b(a+1)}e^{-2\pi i bz}\frac{e^{-3\pi i/4}}{\sqrt{a/b}}\sum_{m=0}^{a-1}e^{-\frac{\pi i b}{a}(z-m+\frac12)^2},
\end{aligned}\end{equation}
\end{proposition}
\begin{proof}
By induction from Proposition \ref{P:Phi} (1) we get
\begin{equation}\label{firsteq}
\Phi(z+a,\tau)-\Phi(z,\tau)=\frac{i}{\sqrt{-i\tau}}\sum_{m=0}^{a-1}e^{\frac{\pi i}{\tau}(z+m+\frac12)^2}.
\end{equation}
From Proposition \ref{P:Phi} (3) we get the next equations, where $q=e^{\pi i \tau}$, next to each equation we put a factor. 
\[\begin{aligned}
&\Phi(z,\tau)+e^{-2\pi i z}q^{-1}\Phi(z+\tau,\tau)=1, &&1;\\
&\Phi(z+\tau,\tau)+e^{-2\pi i z}q^{-3}\Phi(z+2\tau,\tau)=1, &&-e^{-2\pi i z}q^{-1};\\
&\Phi(z+2\tau,\tau)+e^{-2\pi i z}q^{-5}\Phi(z+3\tau,\tau)=1, &&e^{-2\pi i z}q^{-4};\\
&\Phi(z+3\tau,\tau)+e^{-2\pi i z}q^{-7}\Phi(z+4\tau,\tau)=1, &&-e^{-2\pi i z}q^{-9};\\
&\dots\\
&\Phi(z+(b-1)\tau,\tau)+e^{-2\pi i z}q^{-(2b-1)}\Phi(z+b\tau,\tau)=1, &&(-1)^{b-1}e^{-2\pi i z}q^{-(b-1)^2};
\end{aligned}\]
Multiplying each of these equations by the factor to its right and adding, we get 
\begin{equation}\label{secondeq}
\Phi(z,\tau)-(-1)^be^{-2\pi i bz}q^{-b^2}\Phi(z+b\tau)=\sum_{n=0}^{b-1}(-1)^ne^{-2\pi i n z}q^{-n^2}.
\end{equation}
Equations \eqref{firsteq}  and \eqref{secondeq} are valid for $\tau\in\Hip$, but for each fixed $z$ the function $\Phi(z,\tau)$ extends analytically for $\tau$ real not equal $0$. Therefore, we can take $\tau=a/b$, taking $\sqrt{-ia/b}=e^{-\pi i/4}\sqrt{a/b}$.  Then 
$z+b\tau =z+a$. And \eqref{firsteq}, \eqref{secondeq} is a system of equations for $\Phi(z,\tau)$ and $\Phi(z,\tau+a)$. From them we get \eqref{E:integral1}.

To get \eqref{E:integral2}  put $z-a$ instead of $z$ in \eqref{firsteq}. Then put $\tau=-b/a$ and eliminate $\Phi(z-a)$. 
\end{proof}

\begin{remark} If in \eqref{E:integral1} or \eqref{E:integral2} we select $z$ such that the right-hand side vanishes, we get a relation between finite sums. In this way, we may prove the law of quadratic reciprocity. In the next section we will get another application. 
\end{remark}

\section{A finite Fourier expansion}

\begin{lemma}
Let $\tau\in\Q^*$ be a non-null rational number. Let $a$ and $b$ be unique integers such that  $\tau=\frac{a}{b}$ with $b\in\N$ and $a\in\Z^*$ with $a\perp b$. 
\begin{itemize}
\item[(a)] The function $u_\tau\colon\Z\to\C$ defined by $u_\tau(m)=(-1)^m e^{\pi i\tau m^2}$ is periodic with period $2b$ when $a$ is even, and period $b$ for odd $a$. 

\item[(b)] For $a$ even 
\begin{equation}\label{E:evencase}
(-1)^me^{\pi i\tau m^2}=\frac{1}{2\sqrt{b}}\sum_{k\bmod 2b} V_\tau(k) e^{2\pi i \frac{k}{2b}m}.
\end{equation}
\item[(c)] For $a$ odd  the function $u$ has period $b$ and 
\begin{equation}\label{E:oddcase}
(-1)^me^{\pi i\frac{a}{b}m^2}=\frac{1}{\sqrt{b}}\sum_{k\bmod b} V_\tau(k) e^{2\pi i \frac{k}{b}m}.
\end{equation}
\end{itemize}
\end{lemma}

\begin{proof}
In both cases,
\begin{align*}
(-1)^{m+b} e^{\pi i \frac{a}{b}(m+b)^2}&=(-1)^{m} e^{\pi i \frac{a}{b}m^2}\cdot (-1)^b e^{\pi i \frac{a}{b}(2bm+b^2)}\\&=(-1)^{m} e^{\pi i \frac{a}{b}m^2}\cdot (-1)^{b+ab}.
\end{align*}
Hence, for $a$ odd we have $u(m+b)=u(b)$ for even $a$, since in this case $b$ is odd, we have $u(m+b)=-u(m)$ so that $u(m+2b)=u(m)$. 

Equations  \eqref{E:evencase} and \eqref{E:oddcase} are the usual Fourier expansion of a periodic function. By the orthogonality relations, we have 
\begin{equation}\label{E:Ueven}
V_\tau(n)=\frac{1}{2\sqrt{b}}\sum_{m\bmod2b}(-1)^me^{\pi i\frac{a}{b}m^2}e^{-2\pi i\frac{n}{2b}m},\qquad \text{$a$ even},
\end{equation}
and 
\begin{equation}\label{E:Uodd}
V_\tau(n)=\frac{1}{\sqrt{b}}\sum_{m\bmod b}(-1)^me^{\pi i\frac{a}{b}m^2}e^{-2\pi i\frac{n}{b}m},\qquad \text{$a$ odd}.\qedhere
\end{equation}
Therefore, for odd $a$ we have $V_\tau(n)=\widehat{u}_\tau(n)$, but for even $a$ we have 
$V_\tau(n)=\sqrt{2}\,\widehat{u}_\tau(n)$. In this way (see Theorem \ref{root}), we  get $V_\tau(n)$ equal, in all cases, to a root of unity or $0$. 
\end{proof}

\begin{proposition}
Let $\tau\in Q^*$ be such that $\tau=\frac{a}{b}$ with $a\in \Z^*$,  $b\in\N$ and $a\perp b$. Then, for all $n\in\Z$ we have 
\begin{gather}
V_{-\tau}(n)=\overline{V_\tau(n)}, \label{U:-a}\\
V_{\tau+2}(n)=V_\tau(n)\label{U:mod},\\
\begin{cases}
V_\tau(n+b)=V_\tau(n), & \text{for $a$ odd},\\
V_\tau(n+2b)=V_\tau(n), & \text{for $a$ even}.\end{cases}\label{U:periodn}\\
\begin{cases}
V_\tau(b-n)=V_\tau(n),  & \text{for $a$ odd},\\
V_\tau(2b-n)=V_\tau(n),& \text{for $a$ even}.\end{cases}\label{U:final}
\end{gather}
\end{proposition}

\begin{proof}
To proof \eqref{U:-a}, conjugate the equation \eqref{E:Ueven} (or \eqref{E:Uodd}), and then observe that when $m$ runs through representants of the class $\bmod\; 2b$ (or $\bmod\; b$) the numbers $-m$ run also through representants of the same class. The same observation proves \eqref{U:final}

Equations \eqref{U:mod} and \eqref{U:periodn} follow directly from the definition \eqref{E:Ueven} or \eqref{E:Uodd}.
\end{proof}

Quadratic Gaussian sums are defined by 
\[S(a,b)=\sum_{n=0}^{b-1}e^{2\pi i\frac{a}{b}n^2}.\]
We may assume that $a\perp b$, since in the other case, for $d=\gcd(a,b)$, we have 
$S(a,b)=dS(a/d,  b/d)$. When $a\perp b$ the values of $S(a,b)$ can be calculated explicitly; see, for example,  \cite{Hua}*{7.5}. In particular we have 
\begin{proposition}
For $a\perp b$ we have $S(a,b)=0$ when $2\mid b$ but $4\nmid b$. In other cases $(S(a,b)/\sqrt{b})^4=1$ for $b$ odd and $(S(a,b)/2\sqrt{b})^4=1$ for $b$ even.
\end{proposition}

\begin{theorem}\label{root}
Let $\tau\in Q^*$ be such that $\tau=\frac{a}{b}$ with $a\in \Z^*$,  $b\in\N$ and $a\perp b$. When both $a$ and $n$ are even, 
$V_\tau(n)=0$. In all other cases $V_\tau(n)$ is a root of unity. 
\end{theorem}
\begin{proof}
Since $V_{-\tau}(n)=\overline{V_\tau(n)}$, we may assume that $a$ and $b$ are natural numbers.

When $a$ and $n$ are even, the terms in the sum in \eqref{E:Ueven} corresponding to $m$ and $m+b$ are opposite. So, the total sum is $0$. We divide the proof in two cases:

\subsection{Case \texorpdfstring{$a$}{a} odd}
Since $m\mapsto (-1)^me^{\pi i\frac{a}{b}m^2}e^{-2\pi i\frac{n}{b}m}$ is periodic of period $b$ we have 
\[V_\tau(n)=\frac{1}{\sqrt{b}}\sum_{m\bmod b}(-1)^me^{\pi i\frac{a}{b}m^2}e^{-2\pi i\frac{n}{b}m}= \frac{1}{2\sqrt{b}}\sum_{m\bmod 2b}(-1)^me^{\pi i\frac{a}{b}m^2}e^{-2\pi i\frac{n}{b}m}.\]
We have $(-1)^m=(-1)^{m^2}$ and,  therefore,
\[V_\tau(n)=\frac{1}{2\sqrt{b}}\sum_{m\bmod 2b}(-1)^me^{\pi i\frac{a}{b}m^2}e^{-2\pi i\frac{n}{b}m}=\frac{1}{2\sqrt{b}}\sum_{m\bmod 2b}e^{2\pi i\frac{(a+b)m^2-2nm}{2b}}.\]
Since $a\perp b$, there is some $d$ with $(a+b)d\equiv 1\pmod{b}$, then 
\[V_\tau(n)=\frac{1}{2\sqrt{b}}\sum_{m\bmod 2b}e^{2\pi i\frac{(a+b)(m^2-2dnm)}{2b}}=
\frac{e^{-2\pi i (a+b)\frac{d^2n^2}{2b}}}{2\sqrt{b}}\sum_{m\bmod 2b}e^{2\pi i\frac{(a+b)(m-dn)^2}{2b}},\]
$(m-dn)$ run through all the $\bmod 2b$ numbers, so that 
\[V_\tau(n)=\frac{e^{-2\pi i (a+b)\frac{d^2n^2}{2b}}}{2\sqrt{b}}\sum_{m\bmod 2b}e^{2\pi i\frac{(a+b)m^2}{2b}}=\frac{e^{-2\pi i (a+b)\frac{d^2n^2}{2b}}}{2\sqrt{b}} S(a+b,2b).\]
Now, if $b$ is odd $S(a+b,2b)=2S(\frac{a+b}{2},b)$ and this is equal to $2\sqrt{b}$ by a fourth root of unity. If $b$ is even, $S(a+b,2b)$ is equal to $\sqrt{4b} e^{\pi ih/4}$ for some $h\in\Z$. 

\subsection{Case \texorpdfstring{$a$}{a} even.}
In this case $V_\tau(n)=0$ for $n$ even and when $n$  is odd 
\[V_\tau(n)=\frac{1}{2\sqrt{b}}\sum_{m\bmod2b}(-1)^me^{\pi i\frac{a}{b}m^2}e^{-2\pi i\frac{n}{2b}m}=\frac{1}{2\sqrt{b}}\sum_{m\bmod2b}e^{2\pi i \frac{am^2+2\frac{b-n}{2}m}{2b}}.\]
We know that $b$ is odd, so that $\beta=\frac{b-n}{2}$ is an integer. Take $ad\equiv1\bmod b$, then 
\[V_\tau(n)=\frac{1}{2\sqrt{b}}\sum_{m\bmod2b}e^{2\pi i \frac{a(m^2+2d\beta m)}{2b}}=
\frac{e^{-2\pi \frac{a d^2\beta^2}{2b}}}{2\sqrt{b}}\sum_{m\bmod2b}e^{2\pi i \frac{a(m+d\beta)^2}{2b}}=\frac{e^{-2\pi i\frac{a d^2\beta^2}{2b}}}{2\sqrt{b}} S(a,2b).\]
Since $a$ is even, this is 
\[V_\tau(n)=\frac{e^{-2\pi i\frac{a d^2\beta^2}{2b}}}{\sqrt{b}} S(a/2,b).\]
This is also a root of unity. 
\end{proof}

\begin{remark} We may get more explicit values applying the known values for the Gauss sums. Nevertheless, it appears that they can be simplified in some cases. For example,
\[V_{1/b}(n)=(-1)^ne^{\frac{\pi i}{4}}e^{-\frac{\pi i b}{4}}e^{-\frac{\pi i n^2}{b}}.\]
\end{remark}

\section{Riemann auxiliary function with rational parameter}

\begin{theorem}
Let $\tau=-\frac{a}{b}$ be a rational number with  $a\perp b$, with $a$ and $b\in \N$. Then, for $a$ even
\begin{equation}
f(s,\tau)=\frac{1}{\sqrt{b}}\int_{0\swarrow}\frac{e^{-\pi i\tau x^2} x^{-s}}{e^{\pi i x}-e^{-\pi i x}}\Bigl(\sum_{\substack{k\bmod b\\2|k|\le b}}V_\tau(k)\cos\frac{k\pi x}{b}\Bigr)\,dx,
\end{equation}
and for $a$ odd
\begin{equation}
f(s,\tau)=\frac{1}{2\sqrt{b}}\int_{0\swarrow}\frac{e^{-\pi i\tau x^2} x^{-s}}{e^{\pi i x}-e^{-\pi i x}}\Bigl(\sum_{\substack{k\bmod 2b\\|k|\le b}}V_\tau(k)\cos\frac{k\pi x}{2b}\Bigr)\,dx,
\end{equation}
The path of integration  is the line given by $u\mapsto \frac12+ue^{-3\pi i/4}$ (or $u\mapsto \frac12+ue^{\pi i/4}$ when $\tau>0$). 
\end{theorem}

\begin{proof}
It is sufficient to find the expression for $\sum_{m\in\Z}\frac{x e^{\pi i \tau m^2}}{x^2-m^2}$. Since $\tau$ is rational, it is clear that we may express this as a finite combination of trigonometric functions. Consider an odd $a$, then by \eqref{E:oddcase} we have 
\[\frac{1}{2\pi i}\sum_{n\in\Z}\frac{x e^{\pi i \tau n^2}}{x^2-n^2}=\frac{1}{\sqrt{b}}\sum_{k\bmod b}V_\tau(k)\frac{1}{2\pi i}\sum_{m\in\Z}(-1)^m\frac{xe^{2\pi i\frac{k}{b}m}}{x^2-m^2}\]
There is always a selection  of representants $k\bmod b$ such that $2|k|\le b$. Taking this representants we may apply Lemma \ref{L:mittag} which yields 
\[\frac{1}{2\pi i}\sum_{n\in\Z}\frac{x e^{\pi i \tau n^2}}{x^2-n^2}=\frac{1}{\sqrt{b}}\sum_{\substack{k\bmod b\\2|k|\le b}}V_\tau(k)\frac{\cos\frac{2k\pi x}{b}}{e^{\pi i x}-e^{-\pi i x}}.\]
The case $a$ even is treated similarly. 
\end{proof}

\subsection{Examples}

\begin{equation}
\Rzeta_{-1}(s)=\Rzeta(s)=\int_{0\swarrow1}\frac{x^{-s}e^{\pi ix^2}}{e^{\pi i x}-e^{-\pi i x}}\,dx.
\end{equation}

\begin{equation}
\begin{aligned}
\Rzeta_{-4/3}(s)=&\frac{1}{2\sqrt{3}}\int_{0\swarrow1}\frac{x^{-s}e^{\frac{4\pi i}{3}x^2}}{e^{\pi i x}-e^{-\pi i x}}\Bigl(2e^{-\frac{\pi i}{6}}\cos\tfrac{\pi x}{3}+i\cos(\pi x)\Bigr)\,dx\\
&+\frac{1}{4}\int_{0\swarrow1}\frac{x^{-s}e^{\frac{3\pi i}{4}x^2}}{e^{\pi i x}-e^{-\pi i x}}\Bigl(e^{\frac{\pi i}{4}}+2\cos\tfrac{\pi x}{2}+e^{-\frac{3\pi i}{4}}\cos(\pi x)\Bigr)\,dx.
\end{aligned}
\end{equation}

\begin{equation}
\begin{aligned}
\Rzeta_{-3/2}(s)=\frac{1}{2\sqrt{2}}&\int_{0\swarrow1}\frac{x^{-s}e^{\frac{3\pi i}{2}x^2}}{e^{\pi i x}-e^{-\pi i x}}\Bigl(e^{-\frac{\pi i}{4}}+e^{\frac{\pi i}{4}}\cos(\pi x)\Bigr)\,dx\\
&+\frac{1}{2\sqrt{3}}\int_{0\swarrow1}\frac{x^{-s}e^{\frac{2\pi i}{3}x^2}}{e^{\pi i x}-e^{-\pi i x}}\Bigl(2e^{\frac{\pi i}{6}}\cos\frac{\pi x}{3}-i\cos(\pi x)\Bigr)\,dx.
\end{aligned}
\end{equation}

\begin{equation}
\begin{aligned}
\Rzeta_{-2}(s)=\frac12\int_{0\swarrow1}&\frac{x^{-s}e^{2\pi ix^2}\cos\pi x}{e^{\pi i x}-e^{-\pi i x}}\,dx\\&+\frac{1}{2\sqrt{2}}\int_{0\swarrow1}\frac{x^{-s}e^{\frac{\pi i}{2}x^2}}{e^{\pi i x}-e^{-\pi i x}}\Bigl(e^{\frac{\pi i}{4}}+e^{-\frac{\pi i}{4}}\cos(\pi x)\Bigr)\,dx.
\end{aligned}
\end{equation}
\begin{equation}
\begin{aligned}
\Rzeta_{-3}(s)=\frac12&\int_{0\swarrow1}\frac{x^{-s}e^{3\pi ix^2}}{e^{\pi i x}-e^{-\pi i x}}\,dx\\&+\frac{1}{2\sqrt{3}}\int_{0\swarrow1}\frac{x^{-s}e^{\frac{\pi i}{3}x^2}}{e^{\pi i x}-e^{-\pi i x}}\Bigl(e^{\frac{\pi i}{2}}+2e^{-\frac{\pi i}{6}}\cos(\tfrac{2\pi x}{3})\Bigr)\,dx.
\end{aligned}
\end{equation}

\vfill

\begin{figure}[H]
  \includegraphics[width=0.495\hsize]{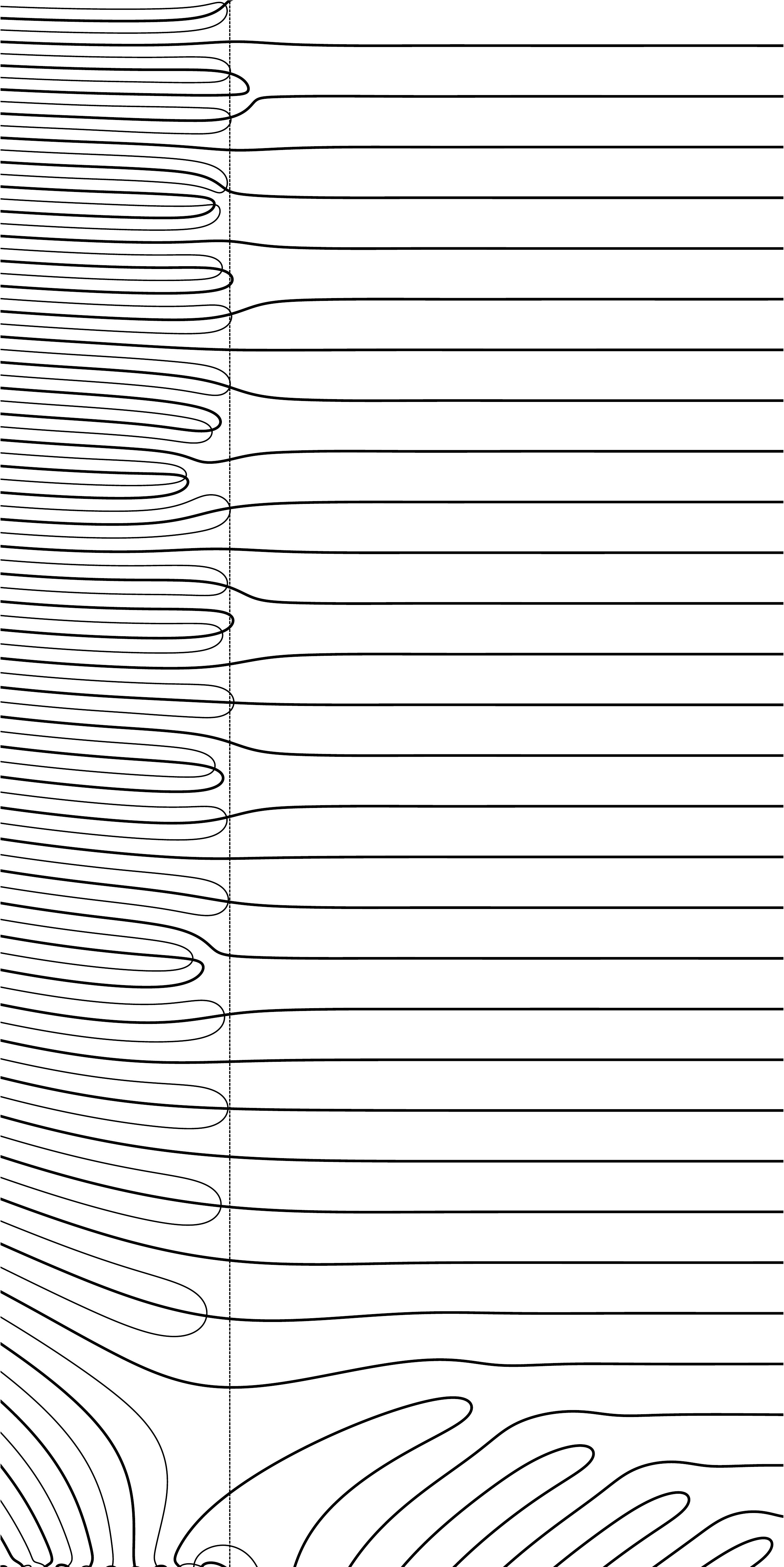}
  \includegraphics[width=0.495\hsize]{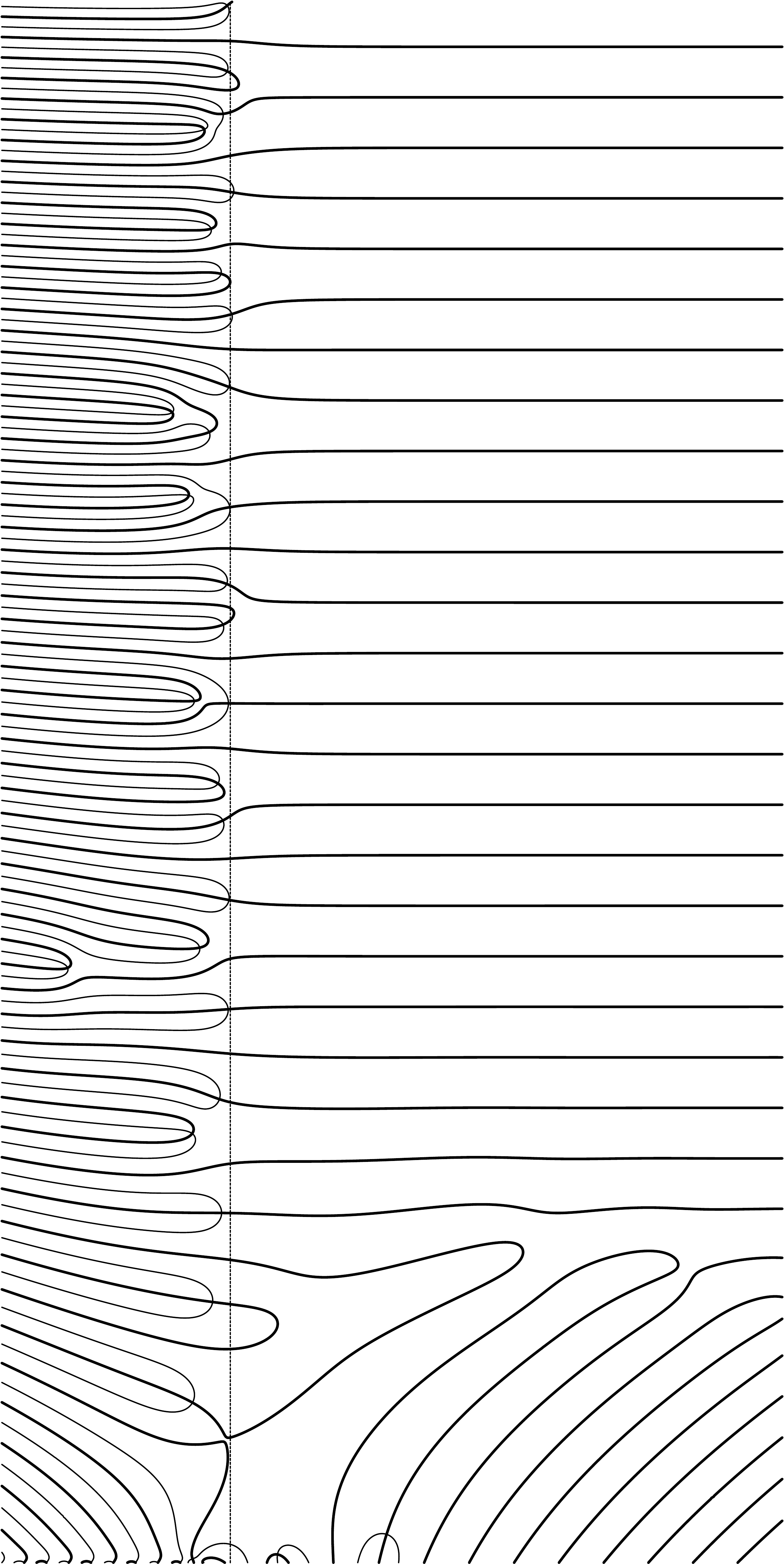}\\
  \caption{$\Rzeta(s)$ and $\Rzeta_{-3/2}(s)$ on $(-20,50)\times(0,140)$}
  \label{F:Xray}
\end{figure}

\vfill

\begin{figure}[H]
  \includegraphics[width=0.495\hsize]{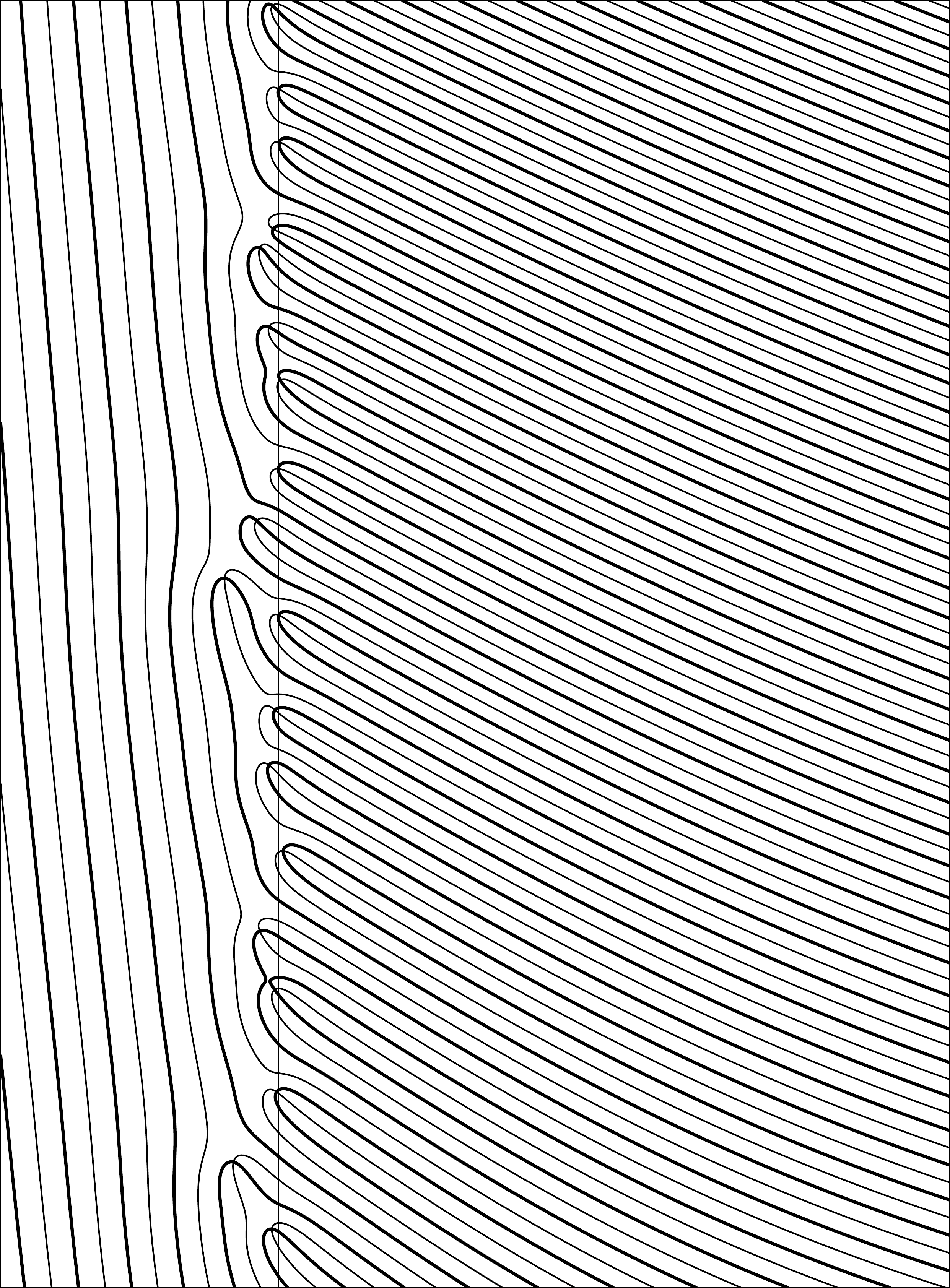}
  \includegraphics[width=0.495\hsize]{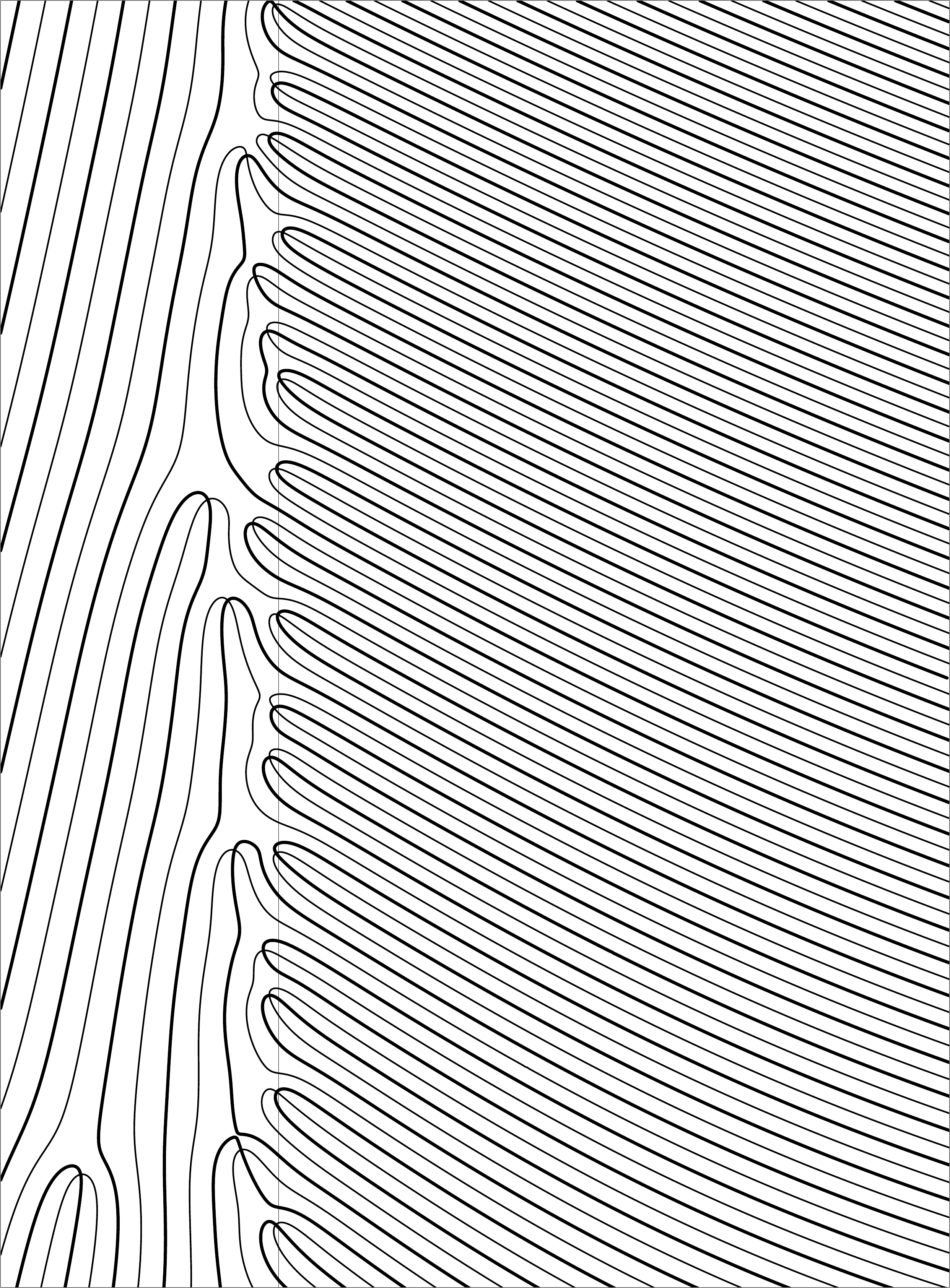}\\
  \includegraphics[width=0.495\hsize]{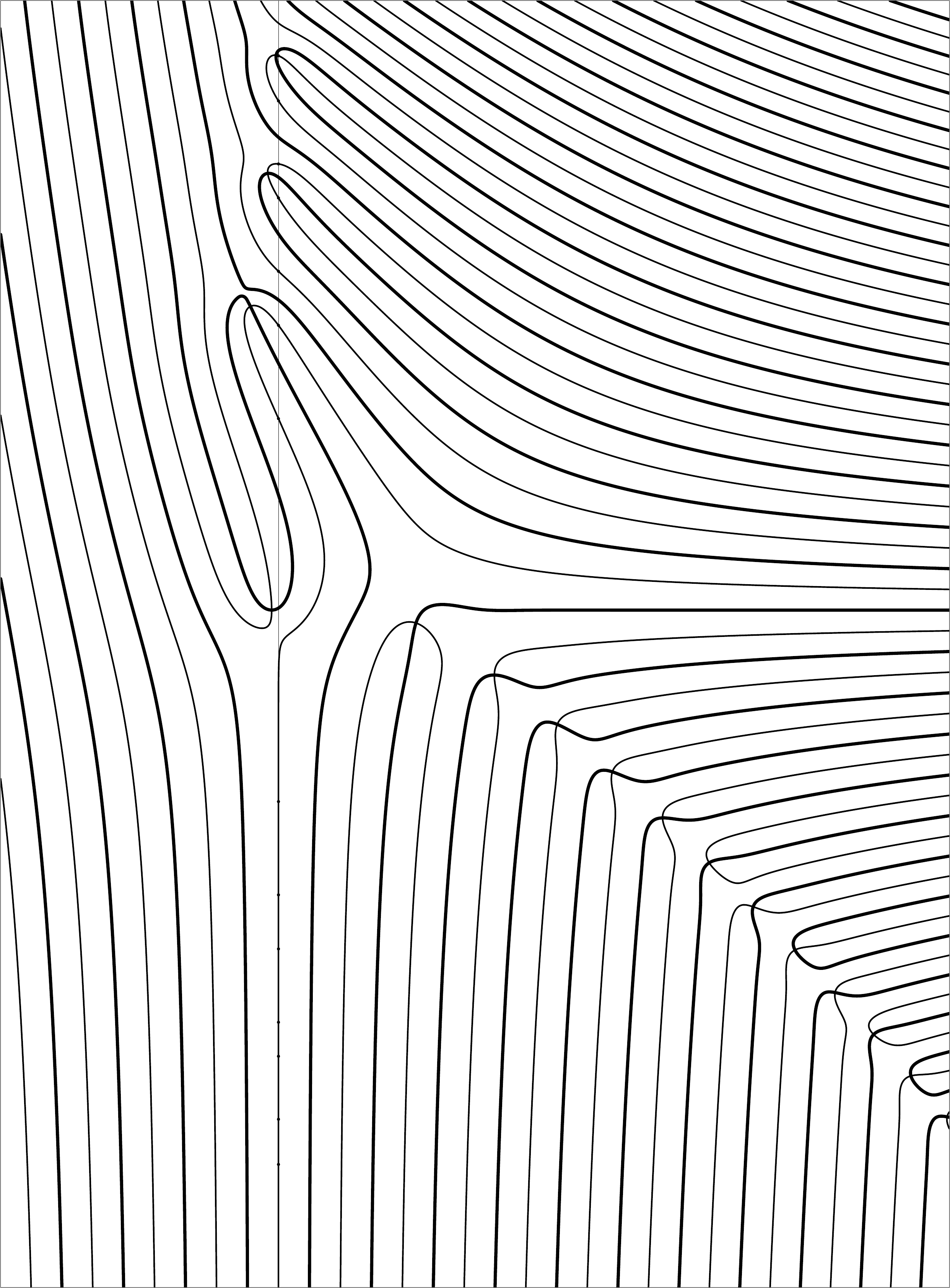}
  \includegraphics[width=0.495\hsize]{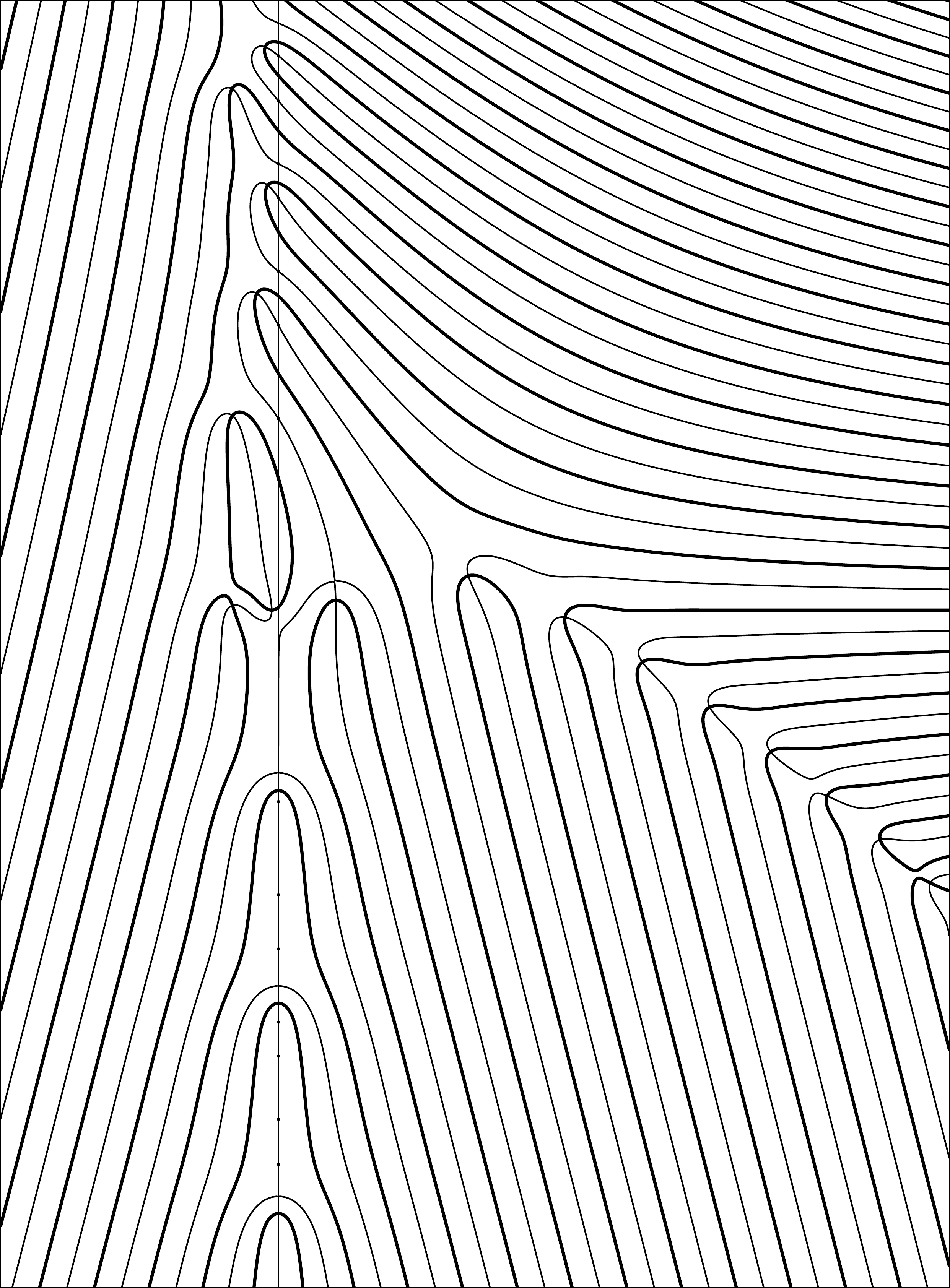}\\
  \caption{$\pi^{-s/2}\Gamma(s/2)\Rzeta(s)$ and $\pi^{-s/2}\Gamma(s/2)\Rzeta_{-3/2}(s)$}
  \label{F:Xray2}
\end{figure}


\noindent Figure \ref{F:Xray2} represents the two functions $\pi^{-\frac{s}{2}}\Gamma(s/2)\Rzeta(s)$ at the left and $\pi^{-\frac{s}{2}}\Gamma(s/2)\Rzeta_{-3/2}(s)$ to the right on $(-20,40)\times(-50,45)$ and  $(-20,40)\times(45,140)$  
\end{document}